%Paper exhibiting p-groups with same integ. coh. groups.
%modified for resubmission on 3/10/94 (originally submitted 10/3/93).  
\magnification \magstep 1
\pageno=1
\def\Bbb #1{{\bf #1}}   %This line needed for printers that don't have 
%                        the AMS fonts for some reason
\def\tilg{{\widetilde G}}
\def\tilh{{\widetilde H}}
\def\bee{B}
\def\co #1 #2{{H^{#1}(\bee{#2})}}
\def\sone{{\Bbb T}}
\def\Hom{{\rm Hom}}
\def\subti #1.{\par\medskip\noindent{\bf #1.}\quad}
\def\proof{\subti Proof.}
\def\arrowright{\rightarrow}
\def\mapright#1{\smash{\mathop{\longrightarrow}\limits^{#1}}}
\def\mapdown#1{\Big\downarrow\rlap{$\vcenter{\hbox{$\scriptstyle#1$}}$}}
\def\qed{$\bullet$\par\medskip}

\def\ca{1}
\def\hu{2}
\def\hv{3}
\def\la{4}
\def\le{5}
\def\stal{6}
\def\stam{7}
\def\ya{8}

\centerline{\bf $p$-groups are not determined by their 
integral cohomology groups}
\medskip
\centerline{Ian Leary}
\centerline{Departement Mathematik,}
\centerline{ETH Zentrum,}
\centerline{8092 Z\"urich.}
\medskip
\centerline{current address: \qquad }
\centerline{Centre de Recerca Matem\`atica,}
\centerline{Institut d'Estudis Catalans,}
\centerline{E 08193 Bellaterra,}
\centerline{Spain}
\medskip

\beginsection Introduction. 

In this paper we prove the result of the title.  More precisely, for each 
prime $p$ we exhibit non-isomorphic $p$-groups $G_1$ and $G_2$ such
that for each $n$ the cohomology groups $H^n(G_1)$ and $H^n(G_2)$ are
isomorphic (throughout this paper the coefficients for cohomology are
the integers unless otherwise stated).  For groups with order
divisible by more than one prime, it has been known for some time that
the integral cohomology ring does not determine the group.  The
following argument is due to Alperin and Atiyah (unpublished).  

If $G$ is an extension with kernel $N$ and quotient $Q$ of coprime
orders, then an easy spectral sequence argument gives
the following ring isomorphism.
$$H^*(G)\cong H^*(N)^Q\otimes H^*(Q)$$
Thus it suffices to find a pair of groups $N$ and $Q$ of coprime
orders such that $Q$ has two non-conjugate actions on $N$ having
isomorphic fixed point subrings in $H^*(N)$.  The smallest such examples
have order 24.  In this case $N$ is cyclic of order three, and $Q$ may be 
any group of order eight having non-isomorphic subgroups $Q_1$ and $Q_2$ 
of order four, where the action of $Q$ on $N$ in the group $G_i$ is given 
by the faithful action of $Q/Q_i$.  
Using a different argument 
Larson has been able to exhibit arbitrarily many metacyclic groups
having isomorphic integral cohomology rings [\la].  

The only previous examples of $p$-groups having isomorphic
integral cohomology groups known to the author are a family of 
(pairs of) $p$-groups for $p$ at least
5 discovered by Yagita [\ya] and a similar family for $p=3$, for which
the author has been able to show further that the groups have
isomorphic integral cohomology rings [\le].  The advantages of this
paper are that a wider class of examples is exhibited, including
$p$-groups for $p=2$, and that very little calculation is involved.
This paper was inspired by my joint work with Nobuaki Yagita, and
relies on a technique suggested to me by Peter Kropholler.  

\beginsection Methods.  

Let $\tilg$ be a compact Lie group whose identity component is central
and isomorphic to $\sone$, the group of complex numbers of unit
modulus.  If $\phi$ is a map from $\tilg$ onto $\sone$, then the
pullback to $\bee\tilg$ of the canonical bundle over $\bee\sone$ is a
principal $\sone$-bundle, with total space the classifying space of
the group $G=\ker(\phi)$.  Under the isomorphism (natural in $\tilg$) 
$\Hom(\tilg,\sone)\cong\co 2 \tilg$, $\phi$ corresponds to the Chern
class of this bundle.  The cohomology of $G$ may be computed from that
of $\tilg$ using the Gysin sequence for
this bundle.  
$$\mapright{}\co {n-2} \tilg\mapright{\times\phi}
\co n \tilg\mapright{}\co n G \mapright{}\co {n-1} \tilg \mapright{}$$
In this long exact sequence, the map from $\co {n-2} \tilg$ to $\co n \tilg$ 
is multiplication by the element of $\co 2 \tilg$ corresponding to $\phi$. 
The idea of using this technique to study the cohomology
of a finite group was suggested by P.~H.~Kropholler and
J.~Huebschmann [\hu,\hv].  

\subti Notation.  From here onwards, $\tilg$ will stand for a Lie
group as above (that is, with central identity component isomorphic to
$\sone$), 
%but also having group of components a finite $p$-group for
%some prime $p$.  
such that the group of components $\tilg/\sone$ is a finite $p$-group 
for some prime $p$.  

Our examples will be pairs of finite subgroups of such a group $\tilg$
corresponding to similar elements in $\Hom(\tilg,\sone)$.  Before
stating our condition for the cohomology groups of two such groups to
be isomorphic in Lemma~4, we require three propositions, one
concerning the cohomology of groups such as $\tilg$, one concerning
the cohomology of their finite subgroups, and another which will allow
us to unfilter a result obtained from a spectral sequence.  

\proclaim Proposition 1.  If $\tilg$ is a compact Lie group as above,
then for $n$ odd $\co n \tilg$ is a finite $p$-group, and for $n$ even
$\co n \tilg$ is the direct sum of a finite $p$-group and an infinite
cyclic group.  If $\theta$ is an element of infinite order in $\co 2
\tilg$, then multiplication by $\theta$ sends any element of infinite
order in $\co n \tilg$ to an element of infinite order in $\co {n+2}
\tilg$.  The exponent of torsion in $\co * \tilg$ is bounded by the
order of the group of components of $\tilg$.  
\par
\proof Let $T$ be the identity component of $\tilg$.  In the Serre
spectral sequence for the fibration 
$$\bee T\arrowright\bee\tilg\arrowright \bee({\tilg/T}),$$
the $E_2$-page is isomorphic to the tensor product of the cohomology
of $T$ and that of $\tilg/T$.  Thus $E_2^{i,j}$ is a finite $p$-group
except for $E_2^{0,2j}$ which is infinite cyclic, and hence the same
is true for $E_n^{i,j}$ for all $n$.  The $E_\infty$-page now yields a
finite filtration of $\co n \tilg$ in which every subquotient is a
finite $p$-group, except the top one, which is infinite cyclic 
if $n$ is even.  An element of
infinite order in $\co {2m} \tilg$ yields a non-zero element of
$E_\infty^{0,2m}$, and the product of such an element with a non-zero
element of $E_\infty^{0,2}$ will be a non-zero element of
$E_\infty^{0,2m+2}$.  Finally, multiplication by the order of the
group of components factors as the composite of the restriction map
from $\tilg$ to $T$ followed by the transfer back to $\tilg$, and $\co
* T$ is torsion free, so this map must annihilate torsion.  
\qed

\subti Remark.  If $G$ is the kernel of a map from $\tilg$ to $\sone$,
then it is easy to see that $G$ is finite if and only if the map is
onto, which occurs if and only if the map has
infinite order when viewed as an element of $\Hom(\tilg,\sone)$.  

\proclaim Propostion 2.  If $G$ is the kernel of a map from $\tilg$
onto $\sone$ and $p^l$ annihilates the torsion in $\co * \tilg$, then 
$p^{2l}\co {2n+1} G$ is trivial, and $p^{2l}\co {2n} G$ is cyclic.  

\proof Consideration of the Gysin sequence for
$\bee G$ as a $\sone$-bundle over $\bee\tilg$ shows that $\co n G$ is 
expressible as an extension with kernel a quotient of $\co n \tilg$ 
and quotient a subgroup of $\co {n-1} \tilg$.  Each of these groups 
must be finite, because $\co n G$ is finite.  By Proposition~1 a finite 
subgroup of $\co {n-1} \tilg$ has exponent dividing $p^l$, and a finite
quotient group of $\co n \tilg$ has exponent $p^l$ except that when $n$ is 
even it may have one cyclic summand of higher order.  The result now follows.
\qed

\proclaim Proposition 3.  If $A$ is a finite abelian $p$-group such
that $p^mA$ is cyclic, then the isomorphism type of $A$ is determined
by the order of $A$ and the orders of $\Bbb Z_{/p^k}\otimes A$ for
$k\leq m$.  

\proof Express $A$ as a sum of cyclic groups, and let $n(k)$ stand for
the number of cyclic summands in this expression of order at least
$p^k$.  By assumption $n(m+1)$ is at most one, and we may recover the
$n(k)$ for $k\leq m$ from the following equation.  
$$p^{n(k)}={|\Bbb Z_{/p^k}\otimes A|/|\Bbb Z_{/p^{k-1}}\otimes A|}$$
This determines all of the cyclic summands of $A$ except the largest
one, which is obtainable now from $|A|$.  
\qed

\proclaim Lemma 4.  Let $\tilg$ be a Lie group as above, and let
$\theta$, $\psi$ be elements of $\Hom(\tilg,\sone)$ such that $\theta$
has infinite order, and $\psi$ has finite order.  Now let
$\phi_1=p^m\theta+\psi$ and $\phi_2=p^m\theta+q\psi$, where $q$ is an
integer coprime to $p$, and let $G_i$ be the kernel of $\phi_i$.  
Then for sufficiently large $m$ the groups $\co n {G_1}$ and $\co n
{G_2}$ are isomorphic for all $n$.  

\proof In fact, if $p^l$ annihilates the torsion in $\co * \tilg$ we
may take any $m$ greater than or equal to $2l$, and for the rest of
the proof we fix some $m$ having this property.  
We shall consider the Gysin sequences for $\bee
G_1$ and $\bee G_2$ as $\sone$-bundles over $\bee\tilg$ for various
choices of coefficients.  First we use the Gysin sequence with 
integer coefficients to show that $\co n {G_1}$ and $\co n {G_2}$ have the
same order.  Let $T^n$ stand for the torsion subgroup of $\co n \tilg$, 
and let $F^n$ be a complement to $T^n$, so that $F^n$ is either infinite 
cyclic or trivial.  If we define an automorphism $f$ of $\co n \tilg$ 
as the sum of the identity map on $F^n$ and multiplication by $q$ on $T^n$, 
then because $p^m$ annihilates $T^n$ the following diagram commutes. 
$$\relax\matrix{\co n \tilg&\mapright{\times\phi_1}&\co {n+2} \tilg\cr
\mapdown{{\rm Id}}&&\mapdown{f}\cr
\co n \tilg&\mapright{\times\phi_2}&\co {n+2} \tilg\cr}$$
It follows that the kernels of $\times\phi_1$ and $\times\phi_2$ have the 
same order, and similarly for the cokernels, and hence that $\co n {G_1}$
and $\co n {G_2}$ have the same order. 

Next consider the Gysin sequences with $\Bbb Z_{/p^k}$ coefficients
for $k\leq m$.  For these sequences the maps from $\co n \tilg$ to 
$\co {n+2} \tilg$ differ only by a unit multiple, and so have identical 
kernels and cokernels.  It follows
\def\zm{{\Bbb Z_{/p^k}}}
that $\co n {G_1;\zm}$ and $\co n {G_2;\zm}$ have the same order for
each $k\leq m$.  Now the universal coefficient theorem (expressing
$\zm$-cohomology in terms of integral cohomology) implies that 
$\zm\otimes\co n {G_1}$ and $\zm \otimes\co n {G_2}$ have the same
order for all $n$ and for all $k\leq m$.  We are now able to apply
Proposition~3 (since $p^m\co n {G_i}$ is cyclic by Propostion~2) and
deduce that for all $n$, $\co n {G_1}$ and $\co n {G_2}$ are
isomorphic.  \qed

\beginsection Examples.  

It remains to construct, for each prime $p$, pairs of non-isomorphic 
$p$-groups satisfying the
conditions of Lemma~4.  First we consider odd primes.  Fix an odd
prime $p$, and $n$ a divisor of $p-1$.  Now for any $m\geq 0$ and any
%$q\in \Bbb F_p^\times$ 
integer $q$ coprime to $p$
define a group $G=G(m,q)$ by the presentation
$$\eqalign{%Change made on 30th Nov.
G=\langle A,B\mid G''=1,\, &A^{p^{m+2}}=B^p=1,\, [A,G']=1,\,\cr
&\qquad[A^p,G]=1,\,[B,\ldots,[B,[B,A]]\ldots]=A^{qp^{m+1}}\rangle,}$$ 
where the letter $B$ occurs $n$ times in the repeated commutator and
$G'$ stands for the derived subgroup of $G$.  The subgroup of $G$
generated by $A$ and the derived subgroup is abelian and isomorphic to
$C_{p^{m+2}}\oplus(C_p)^{n-1}$, and $G$ is the split extension with
kernel this group and quotient of order $p$ generated by $B$. 
It is easy to see that $G(m,q)$ and $G(m,p+q)$ are isomorphic via a
map sending $A$ to $A$ and $B$ to $B$, so that the isomorphism type of 
$G$ depends only on $m$ and the image of $q$ in $\Bbb F_p^\times$. 
Any elements $A'$ and $B'$ generating $G$ such that $A'$ commutes 
with $G'$ and $B'$ has order $p$ must satisfy a similar 
presentation, except that the new $q$ may be the old one multiplied by
%an $n$th power.  
the $n$th power of any integer coprime to $p$.  
It follows that for fixed $m$ there are $n$
isomorphism types of such groups.  For $n=2$ these are the groups
already considered by the author and Yagita [\le,\ya].

For all $m$ and $q$, these groups are normal subgroups of the Lie
group $\tilg$ with presentation 
$$\eqalign{%Change made on 30th Nov. 1994
\tilg = \langle X,Y,\sone \mid \tilg''=1,\, &X^p=Y^p=1,\,
[\sone,\tilg]=1,\cr
&\qquad
[X,\tilg']=1,\, [Y,\ldots,[Y,[Y,X]]\ldots]=\exp(2\pi i/p)\rangle,}$$
and the corresponding homomorphisms from $\tilg$ to $\sone$ are
related as required by Lemma~4.  More precisely, let $\theta$ be the
homomorphism from $\tilg$ to $\sone$ that send $X$ to 1, $Y$ to 1, and
restricts to $\sone$ as $z\mapsto z^p$.  Now let $\psi$ be the
homomorphism from $\tilg$ to $\sone$ that sends $\sone$ to 1, $Y$ to 1
and $X$ to $\exp(-2\pi i/p)$.  For $q$ coprime to $p$, 
the kernel of $p^m\theta+q\psi$ is
generated by $X\exp(2\pi iq/p^{m+2})$ and $Y$, and is isomorphic to 
$G(m,q')$, where $qq'\equiv 1$ modulo $p$.  

Now consider the 2-groups $H(m,q)$ where $m\geq 0$ and $q$ is 1 or 3,
with the following presentations. 
$$H=\langle A,B \mid A^{2^{m+4}}=B^8=1,\, [A^4,H]=1,\, [A,H']=1,\,
[B,[B,A]]=A^{q2^{m+2}}\rangle$$
The groups $H(m,1)$ and $H(m,3)$ are not isomorphic by reasoning
similar to that given above.  Let $\tilh$ be the Lie group with
presentation 
$$\tilh=\langle X,Y,\sone \mid X^4=Y^8=1,\, [\sone,\tilh]=1,\,
[X,\tilh']=1,\, [Y,[Y,X]]=i\rangle.$$
Let $\theta$ from $\tilh$ to $\sone$ be defined by $\theta(X)=1$,
$\theta(Y)=1$, $\theta(z)=z^4$, and let $\psi$ from $\tilh$ to $\sone$
be defined by $\psi(X)=-i$, $\psi(Y)=1$, $\psi(z)=1$.  If $q$ is 1 or
3, then the kernel of 
the map $2^m\theta+q\psi$ is generated by $\exp(\pi iq/2^{m+3})X$ and
$Y$, and is isomorphic to $H(m,q)$.  

\proclaim Corollary 5.  For each prime $p$ there are distinct
$p$-groups with isomorphic integral cohomology groups.  

\proof Apply Lemma 4 to the examples given above.  \qed

\subti Remarks.  In contrast to Corollary~5, a result due to Stallings 
[\stal] and Stammbach [\stam] implies that a map $f:G_1\longrightarrow G_2$
between $p$-groups must be an isomorphism if it induces an isomorphism 
between $\co 1 {G_2;{\Bbb Z/p}}$ and $\co 1 {G_1;{\Bbb Z/p}}$ 
and an injection from $\co 2 {G_2;{\Bbb Z/p}}$ to $\co 2 {G_1;{\Bbb Z/p}}$.
 
The author does not know of a pair of finite groups of
different orders having isomorphic integral cohomology groups.  
In [\ca], Carlson introduced an integer invariant $che(G)$ for a finite 
group $G$.  The invariant can be defined in terms of the additive structure
of $\co * G$, and is a multiple of the order of $G$, but this multiple 
depends on $G$ [\ca].

\beginsection References. 

\frenchspacing

\def\paper#1/#2/#3/#4/#5/(#6) #7--#8/{\item{#1} #2, #3, {\it #4,} {\bf #5}
({\oldstyle#6}) {\oldstyle #7}--{\oldstyle#8}.\par\smallskip}
\def\prepaper#1/#2/#3/#4/#5/{\item{#1} #2, #3, {\it #4} {#5}.
\par\smallskip}

\paper \ca/J. F. Carlson/Exponents of modules and maps/Invent. 
Math./95/(1989) 13--24/

\paper \hu/J.~Huebschmann/Perturbation theory and free resolutions for
nilpotent groups of class 2/J. of Algebra/126/(1989) 348--99/

\paper \hv/J.~Huebschmann/Cohomology of nilpotent groups of class 2/J. of
Algebra/126/(1989) 400--50/

\prepaper \la/D. S. Larson/The integral cohomology rings of split
metacyclic groups/Unpublished report, Univ. of
Minnesota/({\oldstyle1987})/

\prepaper \le/I. J. Leary/3-groups are not determined by their integral 
cohomology rings/J. Pure and Appl. Algebra,/to appear/

\paper \stal/J. Stallings/Homology and central series of groups/J. of
Algebra/2/(1965) 170--81/

\paper \stam/U. Stammbach/Angewandungen der Homologietheorie der
Gruppen auf Zentralreihen und auf Invarianten von Pr\"asentierungen/Math. 
Z./94/(1966) 155--77/

\paper \ya/N. Yagita/Cohomology for groups of $\hbox{rank}_p(G)=2$
and Brown-Peterson cohomology/J. Math. Soc. Japan/45/(1993) 627--44/

%Preprint/({\oldstyle1991})/

\bye